\documentclass[preprints,article,accept,moreauthors,pdftex]{my-mdpi}

\def\R{\mathbb{R}}

\firstpage{1} 
\pubvolume{xx}
\issuenum{1}
\articlenumber{5}
\pubyear{2019}
\copyrightyear{2019}
\history{Submitted 1-Oct-2019; Revised 7-Dec-2019; Accepted 17-Dec-2019 
(\url{https://doi.org/10.3390/mca25010001})}

\pdfoutput=1

\Title{Numerical Optimal Control of HIV Transmission\\ in Octave/MATLAB}

\Author{Carlos Campos $^{1,\dagger}$\orcidA{},
Cristiana J. Silva $^{2,*,\dagger}$\orcidB{} 
and Delfim F. M. Torres $^{2,\dagger}$\orcidC{}}

\AuthorNames{Carlos Campos, Cristiana J. Silva and Delfim F. M. Torres}

\address{%
$^{1}$ \quad Department of Mathematics, Polytechnic of Leiria, 2411-901 Leiria, Portugal; carlos.campos@ipleiria.pt\\
$^{2}$ \quad Center for Research and Development in Mathematics and Applications (CIDMA),\newline 
Department of Mathematics, University of Aveiro, 3810-193 Aveiro, Portugal; delfim@ua.pt
}

\corres{Correspondence: cjoaosilva@ua.pt}

\firstnote{These authors contributed equally to this work.}

\abstract{We provide easy and readable GNU Octave/MATLAB code for the simulation 
of mathematical models described by ordinary differential equations 
and for the solution of optimal control problems through Pontryagin's
maximum principle. For that, we consider a normalized HIV/AIDS transmission 
dynamics model based on the one proposed in our recent contribution 
(Silva, C.J.; Torres, D.F.M. A SICA compartmental model in epidemiology 
with application to HIV/AIDS in Cape Verde. \emph{Ecol. Complex.} 
\textbf{2017}, \emph{30}, 70--75), given by a system of four ordinary 
differential equations. An HIV initial value problem is solved numerically 
using the \texttt{ode45} {GNU} Octave function and three standard methods 
implemented by us in Octave/MATLAB: Euler method and second-order and 
fourth-order Runge--Kutta methods. Afterwards, a control function is introduced 
into the normalized HIV model and an optimal control problem is formulated, 
where the goal is to find the optimal HIV prevention strategy that maximizes 
the fraction of uninfected HIV individuals with the least HIV new infections 
and cost associated with the control measures. The optimal control problem is 
characterized analytically using the Pontryagin Maximum Principle, and 
the extremals are computed numerically by implementing a forward-backward 
fourth-order Runge--Kutta method. Complete algorithms, for both uncontrolled 
initial value and optimal control problems, developed under the free 
GNU Octave software and compatible with MATLAB are provided along the article.}



\begin{document}


\section{Introduction}

In recent years, mathematical modeling of processes in biology and medicine, 
in particular in epidemiology, has led to significant scientific advances 
both in mathematics and biosciences \cite{MR3024808,MR3999697}. 
Applications of mathematics in biology are completely opening new pathways 
of interactions, and this is certainly true in the area of optimal control:
a branch of applied mathematics that deals with finding control laws 
for dynamical systems over a period of time such that an objective 
functional is optimized \cite{Pontryagin,MR3810766}. It has numerous 
applications in both biology and medicine 
\cite{MR2149792,MR3232354,MR3660471,MR3999702}.

To find the best possible control for taking a dynamical system from one state to another, 
one uses, in optimal control theory, the celebrated Pontryagin's maximum principle (PMP), 
formulated in 1956 by the Russian mathematician Lev Pontryagin and his collaborators 
\cite{Pontryagin}. Roughly speaking, PMP states that it is necessary for any optimal control 
along with the optimal state trajectory to satisfy the so-called Hamiltonian system, 
which is a two-point boundary value problem, plus a maximality condition on the Hamiltonian. 
Although a classical result, PMP is usually not taught to biologists
and mathematicians working on mathematical biology.
Here, we show how such scientists can easily implement
the necessary optimality conditions given by the PMP, numerically,
and can benefit from the power of optimal control theory. For that, we 
consider a mathematical model for HIV. 

HIV modeling and optimal control is a subject under strong current research:
see, e.g., Reference \cite{MR3813045} and the references therein. Here, we 
consider the {SICA} epidemic model for HIV transmission proposed 
in References~\cite{EcoComplexity2017,DCDS2018}, formulate
an optimal control problem with the goal to find 
the optimal HIV prevention strategy that maximizes the fraction of 
uninfected HIV individuals with least HIV new infections and cost 
associated with the control measures, and give complete algorithms 
in {GNU} 
Octave to solve the considered problems. We trust that our work,
by providing the algorithms in an open programming language, contributes
to reducing the so-called ``replication crisis'' (an ongoing methodological crisis 
in which it has been found that many scientific studies are difficult 
or impossible to replicate or reproduce~\cite{An2018}) in the area
of optimal biomedical research. We trust our current work will be very useful 
to a practitioner from the disease control area and will become a reference
in the field of epidemiology for those interested to include an optimal 
control component in their work.


\section{A Normalized SICA HIV/AIDS Model}

We consider the SICA epidemic model for HIV transmission proposed 
in References \cite{EcoComplexity2017,DCDS2018}, which is given by the following 
system of ordinary differential equations:
\begin{equation}
\label{eq:model:1}
\begin{cases}
S'(t) = b N(t) - \lambda(t) S(t) - \mu S(t)\\[0.2 cm]
I'(t) =  \lambda(t) S(t) - (\rho + \phi + \mu)I(t) 
+ \alpha A(t)  + \omega C(t) \\[0.2 cm]
C'(t) = \phi I(t) - (\omega + \mu)C(t)\\[0.2 cm]
A'(t) =  \rho \, I(t) - (\alpha + \mu + d) A(t) \, .
\end{cases}
\end{equation}
The model in Equation \eqref{eq:model:1} subdivides human population into four mutually exclusive 
compartments: susceptible individuals ($S$); 
HIV-infected individuals with no clinical symptoms of AIDS 
(the virus is living or developing in the individuals 
but without producing symptoms or only mild ones) 
but able to transmit HIV to other individuals ($I$); 
HIV-infected individuals under {ART} 
treatment (the so-called 
chronic stage) with a viral load remaining low ($C$); 
and HIV-infected individuals with AIDS clinical symptoms ($A$).
The total population at time $t$, denoted by $N(t)$, is given by
$N(t) = S(t) + I(t) + C(t) + A(t)$.
Effective contact with people infected with HIV is at a rate $\lambda(t)$, given by
\begin{equation*}
\lambda(t) = \frac{\beta}{N(t)} \left( I(t) + \eta_C \, C(t)  + \eta_A  A(t) \right),
\end{equation*}
where $\beta$ is the effective contact rate for HIV transmission.
The modification parameter $\eta_A \geq 1$ accounts for the relative
infectiousness of individuals with AIDS symptoms in comparison to those
infected with HIV with no AIDS symptoms (individuals with AIDS symptoms
are more infectious than HIV-infected individuals---pre-AIDS). 
On the other hand, $\eta_C \leq 1$ translates the partial restoration 
of immune function of individuals with HIV infection 
that use ART correctly \cite{EcoComplexity2017}. 
All individuals suffer from natural death at a constant rate $\mu$. 
Both HIV-infected individuals with and without AIDS symptoms have access 
to ART treatment: HIV-infected individuals with no AIDS symptoms $I$ progress to the class 
of individuals with HIV infection under ART treatment $C$ at a rate $\phi$,
and HIV-infected individuals with AIDS symptoms are treated for HIV at rate $\alpha$.
An HIV-infected individual with AIDS symptoms $A$ that starts treatment moves 
to the class of HIV-infected individuals $I$ and after, if the treatment 
is maintained, will be transferred to the chronic class $C$.
Individuals in the class $C$ leave to the class $I$ at a rate $\omega$
due to a default treatment. HIV-infected individuals with 
no AIDS symptoms $I$ that do not take ART treatment progress to the AIDS 
class $A$ at rate $\rho$. Finally, HIV-infected individuals with AIDS 
symptoms $A$ suffer from an AIDS-induced death at a rate $d$. 

In the situation where the total population size $N(t)$ is not constant, 
it is often convenient to consider the proportions 
of each compartment of individuals in the population, namely
$$
s =S/N, \quad i = I/N, \quad  c=C/N, \quad  r =R/N  \, .
$$

The state variables $s$, $i$, $c$, and $a$ satisfy the following system of differential equations:
\begin{equation}
\label{eq:model:2}
\begin{cases}
s'(t) = b(1-s(t)) - \beta(i(t) + \eta_C c(t) + \eta_A a(t)) s(t) + d \, a(t) \, s(t)\\[0.2 cm]
i'(t) =  \beta \left( i(t) + \eta_C \, c(t) + \eta_A a(t) \right) s(t) 
- (\rho + \phi + b )i(t) + \alpha a(t)  + \omega c(t) + d \, a(t)\, i(t) \,  \\[0.2 cm]
c'(t) = \phi i(t) - (\omega + b)c(t) + d \, a(t) \, c(t)\\[0.2 cm]
a'(t) =  \rho \, i(t) - (\alpha + b + d) a(t) + d \, a^2(t)
\end{cases}
\end{equation}
with $s(t) + i(t) + c(t) + a(t) = 1$ for all $t \in [0, T]$.


\section{Numerical Solution of the SICA HIV/AIDS Model}
\label{sec:3}

In this section, we consider Equation \eqref{eq:model:2} 
subject to the initial conditions given by 
\begin{equation}
\label{initcond:num}
 s(0) = 0.6 \, , \quad i(0)=0.2 \, , \quad c(0)=0.1 \, , \quad a(0)=0.1 \, , 
\end{equation}
by the fixed parameter values from Table~\ref{table:parameters}, 
and by the final time value of $T = 20$ (years). 
\begin{table}[H]
\centering
\caption{Parameter values of the HIV/AIDS model in Equation \eqref{eq:model:2} 
taken from Reference \cite{DCDS2018} and references cited therein.}
\label{table:parameters}
\begin{tabular}{l  p{6.5cm} l l} \toprule
{\small{Symbol}} &  {\small{Description}} & {\small{Value}} \\ \midrule
{\small{$\mu$}} & {\small{Natural death rate}} & {\small{$1/69.54$}} \\
{\small{$b$}} & {\small{Recruitment rate}} & {\small{$2.1 \mu$}}  \\
{\small{$\beta$}} & {\small{HIV transmission rate}} & {\small{$1.6$}} \\
{\small{$\eta_C$}} & {\small{Modification parameter}} & {\small{$0.015$}} \\
{\small{$\eta_A$}} & {\small{Modification parameter}} & {\small{$1.3$}} \\	
{\small{$\phi$}} & {\small{HIV treatment rate for $I$ individuals}} &  {\small{$1$}} \\
{\small{$\rho$}} & {\small{Default treatment rate for $I$ individuals}}
& {\small{$0.1 $}} \\
{\small{$\alpha$}} & {\small{AIDS treatment rate}}
& {\small{$0.33 $}} \\
{\small{$\omega$}} & {\small{Default treatment rate for $C$ individuals}}
& {\small{$0.09$}} \\
{\small{$d$}} & {\small{AIDS induced death rate}} & {\small{$1$}} \\ \bottomrule
\end{tabular}
\end{table}

All our algorithms, developed to solve numerically the initial value problems in Equations
\eqref{eq:model:2} and \eqref{initcond:num}, are developed under the free 
GNU Octave software (version 5.1.0), a high-level programming language 
primarily intended for numerical computations and is the major free 
and mostly compatible alternative to MATLAB \cite{octave}.
We implement three standard basic numerical techniques:  
Euler, second-order Runge--Kutta, and fourth-order Runge--Kutta.
We compare the obtained solutions with the one obtained using 
the \texttt{ode45} GNU Octave function. 


\subsection{Default \texttt{ode45} routine of GNU Octave}

Using the provided \texttt{ode45} function of GNU Octave, 
which solves a set of non-stiff ordinary differential equations 
with the well known explicit Dormand--Prince method of order 4, 
one can solve the initial value problems in Equations 
\eqref{eq:model:2} and \eqref{initcond:num} as follows: 


\medskip

{\small

\texttt{function dy = odeHIVsystem(t,y)}

\texttt{\ \% Parameters of the model}

\texttt{\ mi = 1.0 / 69.54; b = 2.1 * mi; beta = 1.6;}

\texttt{\ etaC = 0.015; etaA = 1.3; fi = 1.0; ro = 0.1;}

\texttt{\ alfa = 0.33; omega = 0.09; d = 1.0;}

\texttt{\ }

\texttt{\ \% Differential equations of the model}

\texttt{\ dy = zeros(4,1);}

\texttt{\ aux1 = beta * (y(2) + etaC * y(3) + etaA * y(4)) * y(1);}

\texttt{\ aux2 = d * y(4);}

\texttt{\ dy(1) = b * (1 - y(1)) - aux1 + aux2 * y(1);}

\texttt{\ dy(2) = aux1 - (ro + fi + b - aux2) * y(2) + alfa * y(4) + omega * y(3);}

\texttt{\ dy(3) = fi * y(2) - (omega + b - aux2) * y(3);}

\texttt{\ dy(4) = ro * y(2) - (alfa + b + d - aux2) * y(4);}

}

\medskip

On the GNU Octave interface, one should then type the following instructions:

\medskip

{\small

\texttt{$>$$>$ T = 20; N = 100;}

\texttt{$>$$>$ [vT,vY] = ode45(@odeHIVsystem,[0:T/N:T],[0.6 0.2 0.1 0.1]);}

}

\medskip

Next, we show how such approach compares with standard numerical techniques.


\subsection{Euler's Method}

Given a well-posed initial-value problem
\[
\frac{dy}{dt}=f\left(  t,y\right)  \hspace{0.2in}\text{with}\hspace
{0.2in}y\left(  a\right)  =\alpha\hspace{0.2in}\text{and}\hspace{0.2in}a\leq
t\leq b,
\]
Euler's method constructs a sequence of approximation points $\left(
t,w\right)  \approx\left(  t,y\left(  t\right)  \right)$ to the exact
solution of the ordinary differential equation by
$t_{i+1}=t_{i}+h$ and $w_{i+1}=w_{i}+hf\left(  t_{i},w_{i}\right)$,
$i=0,1,\ldots,N-1$, where $t_{0}=a$, 
$w_{0}=\alpha$, and $h=\left(b-a\right)/N$.
Let us apply Euler's method to approximate each one of the four 
state variables of the system of ordinary
differential equations (Equation \eqref{eq:model:2}). 
Our \texttt{odeEuler} GNU Octave implementation is as follows:


\medskip

{\small

\texttt{function dy = odeEuler(T)}

\texttt{\ \% Parameters of the model}

\texttt{\ mi = 1.0 / 69.54; b = 2.1 * mi; beta = 1.6;}

\texttt{\ etaC = 0.015; etaA = 1.3; fi = 1.0; ro = 0.1;}

\texttt{\ alfa = 0.33; omega = 0.09; d = 1.0;}

\texttt{\ }

\texttt{\ \% Parameters of the Euler method}

\texttt{\ test = -1; deltaError = 0.001; M = 100;}

\texttt{\ t = linspace(0,T,M+1); h = T / M;}

\texttt{\ S = zeros(1,M+1); I = zeros(1,M+1);}

\texttt{\ C = zeros(1,M+1); A = zeros(1,M+1);}

\texttt{\ }

\texttt{\ \% Initial conditions of the model}

\texttt{\ S(1) = 0.6; I(1) = 0.2; C(1) = 0.1; A(1) = 0.1;}

\texttt{\ }

\texttt{\ \% Iterations of the method}

\texttt{\ while(test $<$ 0)}

\texttt{\ \ \ oldS = S; oldI = I; oldC = C; oldA = A;}

\texttt{\ \ \ for i = 1:M}

\texttt{\ \ \ \ \ \ \% Differential equations of the model}

\texttt{\ \ \ \ \ \ aux1 = beta * (I(i) + etaC * C(i) + etaA * A(i)) * S(i);}

\texttt{\ \ \ \ \ \ aux2 = d * A(i);}

\texttt{\ }

\texttt{\ \ \ \ \ \ auxS = b * (1 - S(i)) - aux1 + aux2 * S(i);}

\texttt{\ \ \ \ \ \ auxI = aux1 - (ro + fi + b - aux2) * I(i) + alfa * A(i) + omega * C(i);}

\texttt{\ \ \ \ \ \ auxC = fi * I(i) - (omega + b - aux2) * C(i);}

\texttt{\ \ \ \ \ \ auxA = ro * I(i) - (alfa + b + d - aux2) * A(i);}

\texttt{\ }

\texttt{\ \ \ \ \ \ \% Euler new approximation }

\texttt{\ \ \ \ \ \ S(i+1) = S(i) + h * auxS;}

\texttt{\ \ \ \ \ \ I(i+1) = I(i) + h * auxI;}

\texttt{\ \ \ \ \ \ C(i+1) = C(i) + h * auxC;}

\texttt{\ \ \ \ \ \ A(i+1) = A(i) + h * auxA;}

\texttt{\ \ \ end}

\texttt{\ }

\texttt{\ \ \ \% Absolute error for convergence}

\texttt{\ \ \ temp1 = deltaError * sum(abs(S)) - sum(abs(oldS - S));}

\texttt{\ \ \ temp2 = deltaError * sum(abs(I)) - sum(abs(oldI - I));}

\texttt{\ \ \ temp3 = deltaError * sum(abs(C)) - sum(abs(oldC - C));}

\texttt{\ \ \ temp4 = deltaError * sum(abs(A)) - sum(abs(oldA - A));}

\texttt{\ \ \ test = min(temp1,min(temp2,min(temp3,temp4)));}

\texttt{\ end}

\texttt{\ dy(1,:) = t; dy(2,:) = S; dy(3,:) = I;}

\texttt{\ dy(4,:) = C; dy(5,:) = A;}

}

\medskip

\begin{figure}[H]
\centering
\includegraphics[scale=0.35]{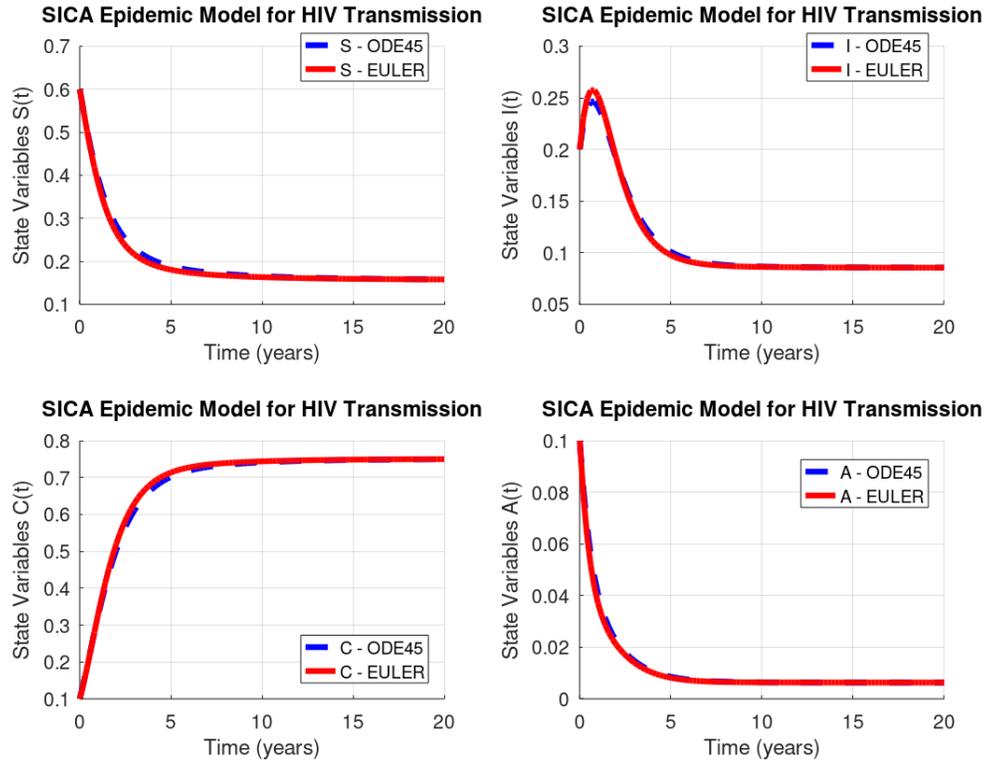}
\caption{HIV/AIDS system (Equation \eqref{eq:model:2}) 
behaviour: {GNU} Octave versus Euler's method.} 
\label{figura1}
\end{figure}

Figure~\ref{figura1}\ shows the solution of the system of ordinary
differential equations (Equation \eqref{eq:model:2}) with the initial 
conditions (Equation \eqref{initcond:num}), 
computed by the \texttt{ode45} GNU Octave function (dashed line) \emph{versus} 
the implemented Euler's method (solid line). As depicted, Euler's method,
although being the simplest method, gives a very good approximation 
to the behaviour of each of the four system variables. Both implementations 
use the same discretization knots in the interval $\left[
0,T\right]$ with a step size given by $h=T/100$.
Euler's method has a global error (total accumulated error) 
of $O\left(  h\right)$, and therefore,
the error bound depends linearly on the step size $h$, which implies that the
error is expected to grow in no worse than a linear manner. Consequently,
diminishing the step size should give correspondingly greater accuracy 
to the approximations. Table~\ref{tabela1}\ lists the norm of the difference vector, 
where each component of this vector is the absolute difference between the results 
obtained by the \texttt{ode45} GNU Octave function 
and our implementation of Euler's method,
calculated by the vector norms $1$, $2$, and $\infty$.
\begin{table}[H]
\centering
\caption{Norms $1$, $2$, and $\infty$ of the difference vector between 
\texttt{ode45} GNU Octave and Euler's results.}
\begin{tabular}
[c]{ccccc}\toprule
System Variables & $S\left(  t\right)  $ & $I\left(  t\right)  $ & $C\left(
t\right)  $ & $A\left(  t\right)  $\\\midrule
$\left\Vert Octave-Euler\right\Vert _{1}$ & $0.4495660$ & $0.1646710$ &
$0.5255950$ & $0.0443340$\\\midrule
$\left\Vert Octave-Euler\right\Vert _{2}$ & $0.0659270$ & $0.0301720$ &
$0.0783920$ & $0.0101360$\\\midrule
$\left\Vert Octave-Euler\right\Vert _{\infty}$ & $0.0161175$ & $0.0113068$ &
$0.0190621$ & $0.0041673$\\\bottomrule
\end{tabular}
\label{tabela1}
\end{table}


\subsection{Runge--Kutta of Order Two}

Given a well-posed initial-value problem, the Runge--Kutta method 
of order two constructs a sequence of approximation points 
$\left(  t,w\right) \approx\left(  t,y\left(  t\right)  \right)$ 
to the exact solution of the ordinary differential equation by
$t_{i+1}=t_{i}+h$, $K_{1}=f\left(  t_{i},w_{i}\right)$,
$K_{2}=f\left(  t_{i+1},w_{i}+hK_{1}\right)$, and
$w_{i+1}=w_{i}+h\dfrac{K_{1}+K_{2}}{2}$,
for each $i=0,1,\ldots,N-1$, where 
$t_{0}=a$, $w_{0}=\alpha$, and $h=\left(b-a\right) /N$.
Our GNU Octave implementation of the Runge--Kutta method of order two applies the
above formulation to approximate each of the four variables 
of the system in Equation \eqref{eq:model:2}. We implement 
the \texttt{odeRungeKutta\_order2} 
function through the following GNU Octave instructions:


\medskip

{\small

\texttt{function dy = odeRungeKutta\_order2(T)}

\texttt{\ \% Parameters of the model}

\texttt{\ mi = 1.0 / 69.54; b = 2.1 * mi; beta = 1.6;}

\texttt{\ etaC = 0.015; etaA = 1.3; fi = 1.0; ro = 0.1;}

\texttt{\ alfa = 0.33; omega = 0.09; d = 1.0;}

\texttt{\ }

\texttt{\ \% Parameters of the Runge-Kutta (2nd order) method}

\texttt{\ test = -1; deltaError = 0.001; M = 100;}

\texttt{\ t = linspace(0,T,M+1); h = T / M; h2 = h / 2;}

\texttt{\ S = zeros(1,M+1); I = zeros(1,M+1);}

\texttt{\ C = zeros(1,M+1); A = zeros(1,M+1);}

\texttt{\ }

\texttt{\ \% Initial conditions of the model}

\texttt{\ S(1) = 0.6; I(1) = 0.2; C(1) = 0.1; A(1) = 0.1;}

\texttt{\ }

\texttt{\ \% Iterations of the method}

\texttt{\ while(test $<$ 0)}

\texttt{\ \ \ oldS = S; oldI = I; oldC = C; oldA = A;}

\texttt{\ \ \ for i = 1:M}

\texttt{\ \ \ \ \ \ \% Differential equations of the model}

\texttt{\ \ \ \ \ \ \% First Runge-Kutta parameter}

\texttt{\ \ \ \ \ \ aux1 = beta * (I(i) + etaC * C(i) + etaA * A(i)) * S(i);}

\texttt{\ \ \ \ \ \ aux2 = d * A(i);}

\texttt{\ \ \ \ \ \ auxS1 = b * (1 - S(i)) - aux1 + aux2 * S(i);}

\texttt{\ \ \ \ \ \ auxI1 = aux1 - (ro + fi + b - aux2) * I(i) + alfa * A(i) + omega * C(i);}

\texttt{\ \ \ \ \ \ auxC1 = fi * I(i) - (omega + b - aux2) * C(i);}

\texttt{\ \ \ \ \ \ auxA1 = ro * I(i) - (alfa + b + d - aux2) * A(i);}

\texttt{\ }

\ \ \\

\texttt{\ \ \ \ \ \ \% Second Runge-Kutta parameter}

\texttt{\ \ \ \ \ \ auxS = S(i) + h * auxS1; auxI = I(i) + h * auxI1;}

\texttt{\ \ \ \ \ \ auxC = C(i) + h * auxC1; auxA = A(i) + h * auxA1;}

\texttt{\ \ \ \ \ \ aux1 = beta * (auxI + etaC * auxC + etaA * auxA) * auxS;}

\texttt{\ \ \ \ \ \ aux2 = d * auxA;}

\texttt{\ \ \ \ \ \ auxS2 = b * (1 - auxS) - aux1 + aux2 * auxS;}

\texttt{\ \ \ \ \ \ auxI2 = aux1 - (ro + fi + b - aux2) * auxI  + alfa * auxA + omega * auxC;}

\texttt{\ \ \ \ \ \ auxC2 = fi * auxI - (omega + b - aux2) * auxC;}

\texttt{\ \ \ \ \ \ auxA2 = ro * auxI - (alfa + b + d - aux2) * auxA;}

\texttt{\ \ \ \ \ \ \% Runge-Kutta new approximation }

\texttt{\ \ \ \ \ \ S(i+1) = S(i) + h2 * (auxS1 + auxS2);}

\texttt{\ \ \ \ \ \ I(i+1) = I(i) + h2 * (auxI1 + auxI2);}

\texttt{\ \ \ \ \ \ C(i+1) = C(i) + h2 * (auxC1 + auxC2);}

\texttt{\ \ \ \ \ \ A(i+1) = A(i) + h2 * (auxA1 + auxA2);}

\texttt{\ \ \ end}

\texttt{\ \ \ \% Absolute error for convergence}

\texttt{\ \ \ temp1 = deltaError * sum(abs(S)) - sum(abs(oldS - S));}

\texttt{\ \ \ temp2 = deltaError * sum(abs(I)) - sum(abs(oldI - I));}

\texttt{\ \ \ temp3 = deltaError * sum(abs(C)) - sum(abs(oldC - C));}

\texttt{\ \ \ temp4 = deltaError * sum(abs(A)) - sum(abs(oldA - A));}

\texttt{\ \ \ test = min(temp1,min(temp2,min(temp3,temp4)));}

\texttt{\ end}

\texttt{\ dy(1,:) = t; dy(2,:) = S; dy(3,:) = I;} \texttt{\ dy(4,:) = C; dy(5,:) = A;}

}

\medskip

\begin{figure}[H]
\centering
\includegraphics[scale=0.30]{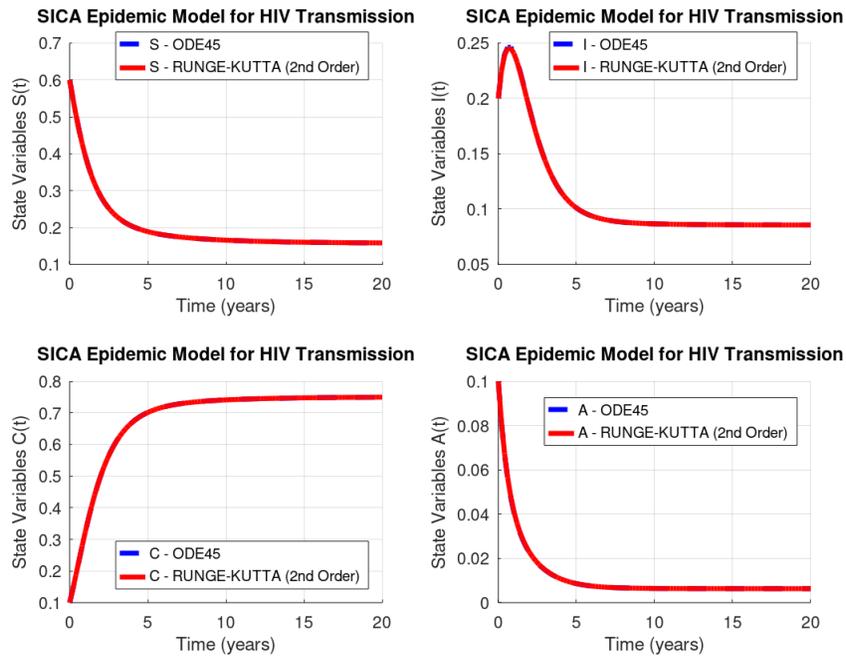}
\caption{HIV/AIDS system (Equation \eqref{eq:model:2}): 
GNU Octave \emph{versus} Runge--Kutta's method of order two.}
\label{figura2}
\end{figure}

Figure~\ref{figura2}\ shows the solution of the system of Equation
\eqref{eq:model:2} with the initial value conditions in Equation \eqref{initcond:num} 
computed by the \texttt{ode45} GNU Octave function (dashed line) 
\emph{versus} our implementation of the Runge--Kutta method
of order two (solid line). As we can see, Runge--Kutta's method produces 
a better approximation than Euler's method, since both curves 
in each plot of Figure~\ref{figura2} are indistinguishable. 
Runge--Kutta's method of order two (RK2) has a global 
truncation error of order $O\left(h^{2}\right)$,
and as it is known, this truncation error at a specified step measures the
amount by which the exact solution to the differential equation fails to
satisfy the difference equation being used for the approximation at that step.
This might seems like an unlikely way to compare the error of various methods,
since we really want to know how well the approximations generated by the
methods satisfy the differential equation, not the other way around. However,
we do not know the exact solution, so we cannot generally determine this, and
the truncation error will serve quite well to determine not only the 
error of a method but also the actual approximation error.
Table~\ref{tabela2} lists the norm of the difference vector
between the results from \texttt{ode45} routine and 
Runge--Kutta's method of order two results.
\begin{table}[!htb]
\centering
\caption{Norms $1$, $2$, and $\infty$ of the difference vector between 
\texttt{ode45} GNU Octave and RK2 results.}
\begin{tabular}
[c]{ccccc}\toprule
System Variables & $S\left(  t\right)  $ & $I\left(  t\right)  $ & $C\left(
t\right)  $ & $A\left(  t\right)  $\\\midrule
$\left\Vert Octave-RungeKutta2\right\Vert _{1}$ & $0.0106530$ & $0.0105505$ &
$0.0151705$ & $0.0044304$\\\midrule
$\left\Vert Octave-RungeKutta2\right\Vert _{2}$ & $0.0014868$ & $0.0025288$ &
$0.0022508$ & $0.0011695$\\\midrule
$\left\Vert Octave-RungeKutta2\right\Vert _{\infty}$ & $0.0003341$ &
$0.0009613$ & $0.0006695$ & $0.0004678$\\\bottomrule
\end{tabular}
\label{tabela2}
\end{table}


\subsection{Runge--Kutta of Order Four}

The Runge--Kutta method of order four (RK4) constructs 
a sequence of approximation points $\left(  t,w\right)
\approx\left(  t,y\left(  t\right)  \right)$ to the exact solution 
of an ordinary differential equation by
$t_{i+1}=t_{i}+h$, $K_{1}=f\left(  t_{i},w_{i}\right)$,
$K_{2}=f\left(  t_{i}+\frac{h}{2},w_{i}+\frac{h}{2}K_{1}\right)$, 
$K_{3}=f\left(  t_{i}+\frac{h}{2},w_{i}+\frac{h}{2}K_{2}\right)$,
$K_{4}=f\left(  t_{i+1},w_{i}+hK_{3}\right)$, and
$w_{i+1}=w_{i}+\dfrac{h}{6}(K_{1}+2K_{2}+2K_{3}+K_{4})$,
for each $i=0,1,\ldots,N-1$, where $t_{0}=a$, 
$w_{0}=\alpha$, and $h=\left(b-a\right)/N$.
Our GNU Octave implementation of the Runge--Kutta method of order four 
applies the above formulation to approximate the solution of  
the system in Equation \eqref{eq:model:2} with the initial conditions 
of Equation \eqref{initcond:num} through the following instructions:


\medskip

{\small

\texttt{function dy = odeRungeKutta\_order4(T)}

\texttt{\ \% Parameters of the model}

\texttt{\ mi = 1.0 / 69.54; b = 2.1 * mi; beta = 1.6;}

\texttt{\ etaC = 0.015; etaA = 1.3; fi = 1.0; ro = 0.1;}

\texttt{\ alfa = 0.33; omega = 0.09; d = 1.0;}

\texttt{\ }

\texttt{\ \% Parameters of the Runge-Kutta (4th order) method}

\texttt{\ test = -1; deltaError = 0.001; M = 100;}

\texttt{\ t = linspace(0,T,M+1);}

\texttt{\ h = T / M; h2 = h / 2; h6 = h / 6;}

\texttt{\ S = zeros(1,M+1); I = zeros(1,M+1);}

\texttt{\ C = zeros(1,M+1); A = zeros(1,M+1);}

\texttt{\ }

\texttt{\ \% Initial conditions of the model}

\texttt{\ S(1) = 0.6; I(1) = 0.2; C(1) = 0.1; A(1) = 0.1;}

\texttt{\ \% Iterations of the method}

\texttt{\ while(test $<$ 0)}

\texttt{\ \ \ oldS = S; oldI = I; oldC = C; oldA = A;}

\texttt{\ \ \ for i = 1:M}

\texttt{\ \ \ \ \ \ \% Differential equations of the model}

\texttt{\ \ \ \ \ \ \% First Runge-Kutta parameter}

\texttt{\ \ \ \ \ \ aux1 = beta * (I(i) + etaC * C(i) + etaA * A(i)) * S(i);}

\texttt{\ \ \ \ \ \ aux2 = d * A(i);}

\texttt{\ \ \ \ \ \ auxS1 = b * (1 - S(i)) - aux1 + aux2 * S(i);}

\texttt{\ \ \ \ \ \ auxI1 = aux1 - (ro + fi + b - aux2) * I(i) + alfa * A(i) + omega * C(i);}

\texttt{\ \ \ \ \ \ auxC1 = fi * I(i) - (omega + b - aux2) * C(i);}

\texttt{\ \ \ \ \ \ auxA1 = ro * I(i) - (alfa + b + d - aux2) * A(i);}

\texttt{\ }

\texttt{\ \ \ \ \ \ \% Second Runge-Kutta parameter}

\texttt{\ \ \ \ \ \ auxS = S(i) + h2 * auxS1; auxI = I(i) + h2 * auxI1;}

\texttt{\ \ \ \ \ \ auxC = C(i) + h2 * auxC1; auxA = A(i) + h2 * auxA1;}

\texttt{\ \ \ \ \ \ aux1 = beta * (auxI + etaC * auxC + etaA * auxA) * auxS;}

\texttt{\ \ \ \ \ \ aux2 = d * auxA;}

\texttt{\ }

\texttt{\ \ \ \ \ \ auxS2 = b * (1 - auxS) - aux1 + aux2 * auxS;}

\texttt{\ \ \ \ \ \ auxI2 = aux1 - (ro + fi + b - aux2) * auxI + alfa * auxA + omega * auxC;}

\texttt{\ \ \ \ \ \ auxC2 = fi * auxI - (omega + b - aux2) * auxC;}

\texttt{\ \ \ \ \ \ auxA2 = ro * auxI - (alfa + b + d - aux2) * auxA;}

\texttt{\ }

\texttt{\ \ \ \ \ \ \% Fird Runge-Kutta parameter}

\texttt{\ \ \ \ \ \ auxS = S(i) + h2 * auxS2; auxI = I(i) + h2 * auxI2;}

\texttt{\ \ \ \ \ \ auxC = C(i) + h2 * auxC2; auxA = A(i) + h2 * auxA2;}

\texttt{\ \ \ \ \ \ aux1 = beta * (auxI + etaC * auxC + etaA * auxA) * auxS;}

\texttt{\ \ \ \ \ \ aux2 = d * auxA;}

\texttt{\ }

\texttt{\ \ \ \ \ \ auxS3 = b * (1 - auxS) - aux1 + aux2 * auxS;}

\texttt{\ \ \ \ \ \ auxI3 = aux1 - (ro + fi + b - aux2) * auxI + alfa * auxA + omega * auxC;}

\texttt{\ \ \ \ \ \ auxC3 = fi * auxI - (omega + b - aux2) * auxC;}

\texttt{\ \ \ \ \ \ auxA3 = ro * auxI - (alfa + b + d - aux2) * auxA;}

\texttt{\ }

\texttt{\ \ \ \ \ \ \% Fourth Runge-Kutta parameter}

\texttt{\ \ \ \ \ \ auxS = S(i) + h * auxS3; auxI = I(i) + h * auxI3;}

\texttt{\ \ \ \ \ \ auxC = C(i) + h * auxC3; auxA = A(i) + h * auxA3;}

\texttt{\ \ \ \ \ \ aux1 = beta * (auxI + etaC * auxC + etaA * auxA) * auxS;}

\texttt{\ \ \ \ \ \ aux2 = d * auxA;}

\texttt{\ }

\texttt{\ \ \ \ \ \ auxS4 = b * (1 - auxS) - aux1 + aux2 * auxS;}

\texttt{\ \ \ \ \ \ auxI4 = aux1 - (ro + fi + b - aux2) * auxI + alfa * auxA + omega * auxC;}

\texttt{\ \ \ \ \ \ auxC4 = fi * auxI - (omega + b - aux2) * auxC;}

\texttt{\ \ \ \ \ \ auxA4 = ro * auxI - (alfa + b + d - aux2) * auxA;}

\texttt{\ }

\texttt{\ \ \ \ \ \ \% Runge-Kutta new approximation }

\texttt{\ \ \ \ \ \ S(i+1) = S(i) + h6 * (auxS1 + 2 * (auxS2 + auxS3) +
	auxS4);}

\texttt{\ \ \ \ \ \ I(i+1) = I(i) + h6 * (auxI1 + 2 * (auxI2 + auxI3) +
	auxI4);}

\texttt{\ \ \ \ \ \ C(i+1) = C(i) + h6 * (auxC1 + 2 * (auxC2 + auxC3) +
	auxC4);}

\texttt{\ \ \ \ \ \ A(i+1) = A(i) + h6 * (auxA1 + 2 * (auxA2 + auxA3) +
	auxA4);}

\texttt{\ \ \ end}

\texttt{\ }

\texttt{\ \ \ \% Absolute error for convergence}

\texttt{\ \ \ temp1 = deltaError * sum(abs(S)) - sum(abs(oldS - S));}

\texttt{\ \ \ temp2 = deltaError * sum(abs(I)) - sum(abs(oldI - I));}

\texttt{\ \ \ temp3 = deltaError * sum(abs(C)) - sum(abs(oldC - C));}

\texttt{\ \ \ temp4 = deltaError * sum(abs(A)) - sum(abs(oldA - A));}

\texttt{\ \ \ test = min(temp1,min(temp2,min(temp3,temp4)));}

\texttt{\ end}

\texttt{\ dy(1,:) = t; dy(2,:) = S; dy(3,:) = I;}

\texttt{\ dy(4,:) = C; dy(5,:) = A;}

}

\medskip

\begin{figure}[H]
\centering
\includegraphics[scale=0.38]{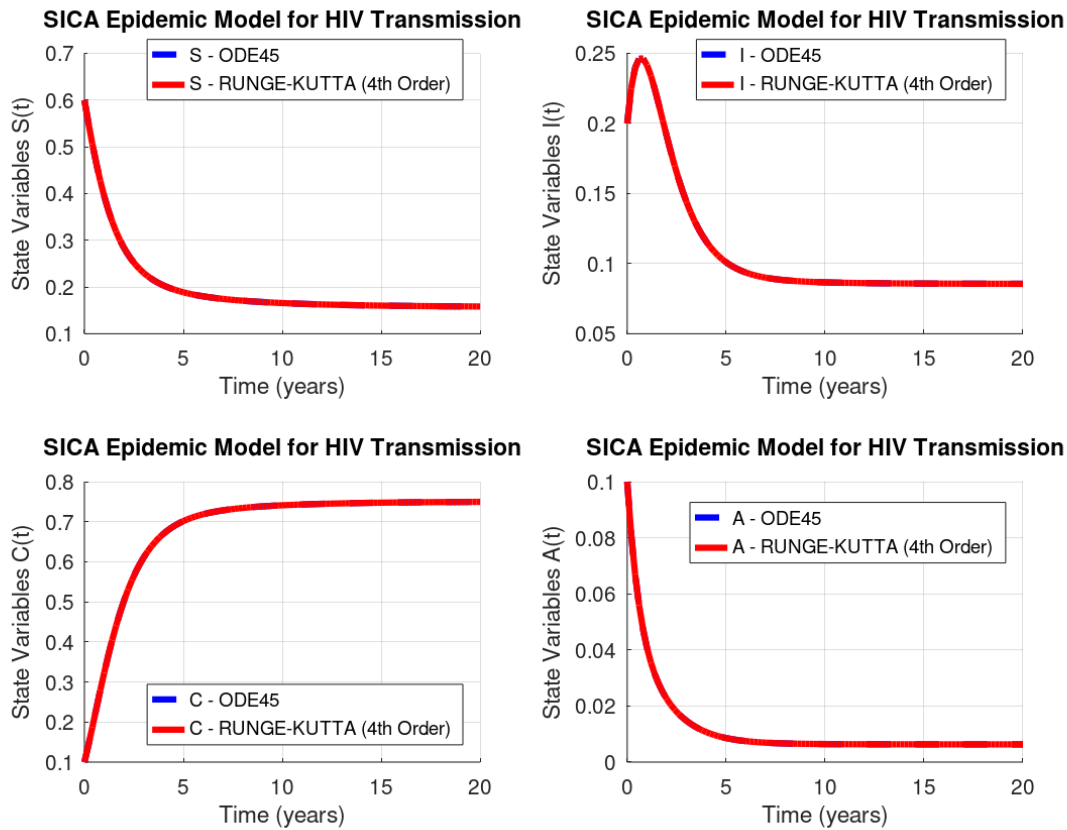}
\caption{HIV/AIDS system (Equation \eqref{eq:model:2}): 
GNU Octave versus Runge--Kutta's method of order four.}
\label{figura3}
\end{figure}

Figure~\ref{figura3}\ shows the solution of the initial value problem in Equations 
\eqref{eq:model:2} and \eqref{initcond:num} computed by the \texttt{ode45} GNU
Octave function (dashed line) \emph{versus} our implementation of 
the Runge--Kutta method of order four (solid line). The results 
of the Runge--Kutta method of order four are extremely good. Moreover, 
this method requires four evaluations per step and
its global truncation error is $O\left(  h^{4}\right)$. 
Table~\ref{tabela3}\ lists the norm of the difference vector
between results obtained by the Octave routine \texttt{ode45} 
and the 4th order Runge--Kutta method.
\begin{table}[H]
\centering
\caption{Norms $1$, $2$, and $\infty$ of the difference vector between 
	\texttt{ode45} GNU Octave and RK4 results.}
\begin{tabular}
[c]{ccccc}\toprule
System Variables & $S\left(  t\right)  $ & $I\left(  t\right)  $ & $C\left(
t\right)  $ & $A\left(  t\right)  $\\\midrule
$\left\Vert Octave-RungeKutta4\right\Vert _{1}$ & $0.0003193$ & $0.0002733$ &
$0.0004841$ & $0.0000579$\\\midrule
$\left\Vert Octave-RungeKutta4\right\Vert _{2}$ & $0.0000409$ & $0.0000395$ &
$0.0000674$ & $0.0000098$\\\midrule
$\left\Vert Octave-RungeKutta4\right\Vert _{\infty}$ & $0.0000107$ &
$0.0000140$ & $0.0000186$ & $0.0000042$\\\bottomrule
\end{tabular}
\label{tabela3}
\end{table}


\section{Optimal Control of HIV Transmission}

In this section, we propose an optimal control problem that will be solved 
numerically in Octave/MATLAB in Section~\ref{sec:oc:numer}. 
We introduce a control function  $u(\cdot)$ in the model 
of Equation \eqref{eq:model:2}, which represents the effort 
on HIV prevention measures, such as condom use 
(used consistently and correctly during every sex act) or oral 
pre-exposure prophylasis (PrEP). The control system is given by
\begin{equation}
\label{eq:model:cont}
\begin{cases}
s'(t) = b(1-s(t)) - (1-u(t))\beta(i(t) 
+ \eta_C c(t) + \eta_A a(t)) s(t) + d \, a(t) \, s(t)\\[0.2 cm]
i'(t) =  (1-u(t))\beta \left( i(t) + \eta_C \, c(t) 
+ \eta_A a(t) \right) s(t) - (\rho + \phi + b )i(t) 
+ \alpha a(t)  + \omega c(t) + d \, a(t)\, i(t) \,  \\[0.2 cm]
c'(t) = \phi i(t) - (\omega + b)c(t) + d \, a(t) \, c(t)\\[0.2 cm]
a'(t) =  \rho \, i(t) - (\alpha + b + d) a(t) + d \, a^2(t),
\end{cases}
\end{equation}
where the control $u(\cdot)$ is bounded between 0 and $u_{\max}$, with $u_{\max} < 1$. 
When the control vanishes, no extra preventive measure for HIV transmission 
is being used by susceptible individuals. We assume that $u_{max}$ is never 
equal to 1, since it makes the model more realistic from a medical point of view.  

The goal is to find the optimal value $u^*$ of the control $u$
along time, such that the associated state trajectories $s^*$, $i^*$, $c^*$, and $a^*$
are solutions of the system in Equation \eqref{eq:model:cont} in the time interval $[0, T]$ with
the following initial given conditions:
\begin{equation}
\label{eq:init}
s(0) \geq 0 \, , \quad i(0) \geq 0 \, , 
\quad c(0) \geq 0 \, , \quad a(0) \geq 0 \, ,
\end{equation}
and $u^*(\cdot)$ maximizes the objective functional given by
\begin{equation}
\label{eq:cost}
J(u(\cdot)) =  \int_0^T \left(s(t) - i(t) - u^2(t)\right) \, dt,
\end{equation}
which considers the fraction of susceptible individuals ($s$) 
and HIV-infected individuals without AIDS symptoms ($i$) and 
the cost associated with the support of HIV transmission measures ($u$). 

The control system in Equation \eqref{eq:model:cont} of ordinary differential 
equations in $\R^4$ is considered with the set of admissible control 
functions given by
\begin{equation}
\label{eq:admiss:cont}
\Omega = \{ u(\cdot) \in L^\infty(0, T) \, | \, 0 \leq u(t) 
\leq u_{max} \, , \, \forall t \in [0, T]  \} \, .
\end{equation}
We consider the optimal control problem of determining 
$\left(s^*(\cdot), i^*(\cdot), c^*(\cdot), a^*(\cdot)\right)$ associated 
to an admissible control $u^*(\cdot) \in \Omega$ on the time interval $[0, T]$, 
satisfying Equation \eqref{eq:model:cont} and the initial conditions 
of Equation \eqref{eq:init} and maximizing the cost functional of Equation \eqref{eq:cost}:
\begin{equation}
\label{eq:cost:max}
J(u^*(\cdot)) = \max_{\Omega} J(u(\cdot)) \, .
\end{equation}
Note that we are considering a $L^2$-cost function: the integrand of the cost functional 
$J$ is concave with respect to the control $u$. Moreover, the control 
system of Equation \eqref{eq:model:cont} is Lipschitz with respect to the state variables $(s, i, c, a)$.
These properties ensure the existence of an optimal control $u^*(\cdot)$ 
of the optimal control problem in Equations \eqref{eq:model:cont}--\eqref{eq:cost:max} 
(see, e.g., Reference \cite{Cesari}).  

To solve optimal control problems, 
two approaches are possible: direct and indirect.
Direct methods consist in the discretization 
of the optimal control problem,
reducing it to a nonlinear programming problem \cite{MR3854267,MR3955185}.
For such an approach, one only needs to use the Octave/MATLAB 
\texttt{fmincon} routine. Indirect methods are more sound
because they are based on Pontryagin's Maximum Principle
but less widespread since they are not immediately available 
in Octave/MATLAB. Here, we show how one can use
Octave/MATLAB to solve optimal control problems
through Pontryagin's Maximum Principle,
reducing the optimal control problem to the solution of 
a boundary value problem. 


\subsection{Pontryagin's Maximum Principle}

According to celebrated Pontryagin's Maximum Principle (see, e.g., Reference \cite{Pontryagin}), 
if $u^*(\cdot)$ is optimal for Equations \eqref{eq:model:cont}--\eqref{eq:cost:max} 
with fixed final time $T$, then there exists a nontrivial absolutely 
continuous mapping $\Lambda : [0, T] \to \R^4$, $\Lambda(t) 
= \left( \lambda_1(t), \lambda_2(t), \lambda_3(t), \lambda_4(t) \right)$, 
called the \emph{adjoint vector}, such that
$$
s' = \frac{\partial H}{\partial \lambda_1}, \ \  
i' = \frac{\partial H}{\partial \lambda_2}, \ \ 
c' = \frac{\partial H}{\partial \lambda_3}, \ \ 
a' = \frac{\partial H}{\partial \lambda_4},  \ \
\lambda_1' =- \frac{\partial H}{\partial s}, \ \ 
\lambda_2' =- \frac{\partial H}{\partial i}, \ \ 
\lambda_3'=- \frac{\partial H}{\partial c}, \ \ 
\lambda_4' =- \frac{\partial H}{\partial a}, 
$$
where 
\begin{equation*}
\begin{split}
H &= H\left(s(t), i(t), c(t), a(t), \Lambda(t), u(t)\right) 
= s(t) - i(t) - u^2(t)\\ 
&+ \lambda_1(t) \Bigl( b(1-s(t)) - (1-u(t))\beta(i(t) + \eta_C c(t) 
+ \eta_A a(t)) s(t) + d \, a(t) \, s(t)  \Bigr)\\
&+ \lambda_2(t) \Bigl(  (1-u(t))\beta \left( i(t) + \eta_C \, c(t) 
+ \eta_A a(t) \right) s(t) - (\rho + \phi + b )i(t) 
+ \alpha a(t)  + \omega c(t) + d \, a(t)\, i(t) \Bigr)\\
&+ \lambda_3(t) \Bigl( \phi i(t) - (\omega + b)c(t) + d \, a(t) \, c(t)  \Bigr)\\
&+ \lambda_4(t) \Bigl( \rho \, i(t) - (\alpha + b + d) a(t) + d \, a^2(t) \Bigr) 
\end{split} 
\end{equation*}
is called the \emph{Hamiltonian} and the maximality condition 
$$
H(s^*(t), i^*(t), c^*(t), a^*(t), \Lambda(t), u^*(t)) 
= \max_{ 0 \leq u \leq u_{max}} H(s^*(t), i^*(t), c^*(t), a^*(t), \Lambda(t), u) 
$$
holds almost everywhere on $[0, T]$. Moreover, the transversality conditions
$$
\lambda_i(T) = 0 \, , \quad \quad i = 1, \ldots , 4,
$$
hold. Applying the Pontryagin maximum principle to the optimal control problem in Equations 
\eqref{eq:model:cont}--\eqref{eq:cost:max}, the following theorem follows. 

\begin{Theorem}
\label{theo:PMP}
The optimal control problem of Equations \eqref{eq:model:cont}--\eqref{eq:cost:max} with fixed final time $T$
admits a unique optimal solution $\left(s^*(\cdot), i^*(\cdot), c^*(\cdot), a^*(\cdot) \right)$ 
associated to the optimal control $u^*(\cdot)$ on $[0, T]$ described by
\begin{equation}
\label{eq:opt:cont}
u^*(t)  =\min\left\{  \max\left\{  0,\frac{\beta\left(
i^*(t)+\eta_{C} c^*(t) + \eta_{A} a^*(t)\right) 
s^*(t) \left(  \lambda_{1}(t)-\lambda_{2}(t) \right)  }{2}\right\}, u_{\max} \right\} \, ,
\end{equation}
where the adjoint functions satisfy 
\begin{eqnarray}
\label{eq:adjointsys}
\begin{cases}
\lambda'_1(t) = -1+\lambda_{1}(t)  \left( b+\left(  1-u^*(t)\right)  \beta\left( i^*(t)
+\eta_{C} c^*(t) +\eta_{A}a^*(t) \right)  - d a^*(t) \right),\\  
\qquad\qquad -\lambda_{2}(t) \left( 1-u^*(t)\right)  \beta\left( i^*(t)+\eta_{C}c^*(t)+\eta_{A}a^*(t)\right)\\
\lambda'_2(t) = 1+\lambda_{1}^*(t)  \left(  1-u^*(t)\right)  \beta s^*(t)
-\lambda_{2}(t) \left( \left(  1-u^*(t)\right)  \beta s^*(t)-\left(  \rho+\phi+s^*(t)\right)
+d a^*(t)\right)  \\
\qquad\qquad -\lambda_{3}(t) \phi-\lambda_{4}(t) \rho,\\
\lambda'_3(t) = \lambda_{1}(t) \left(  1-u^*(t)\right)  \beta\eta_{C}s^*(t)
-\lambda_{2}(t) \left( \left( 1-u^*(t)\right)  \beta\eta_{C}s^*(t)+\omega\right)
+\lambda_{3}(t) \left( \omega+b-d a^*(t)\right), \\
\lambda'_4(t) = \lambda_{1}(t) \left( \left(1-u^*(t)\right)  \beta\eta_{A}s^*(t) 
+ d s^*(t)\right)-\lambda_{2}(t) \left(  \left( 1-u^*(t)\right) \beta \eta_{A} s^*(t) 
+ \alpha + d i^*(t)\right) \\
\qquad\qquad -\lambda_{3}(t) d c^*(t) 
+ \lambda_{4}(t) \left( \alpha + b + d - 2 d a^*(t)\right),
\end{cases}
\end{eqnarray}
subject to the transversality conditions $\lambda_i(T) = 0$, $i = 1, \ldots , 4$. 
\end{Theorem}

\begin{Remark}
The uniqueness of the optimal control $u^*$ is due to the boundedness of the state and adjoint
functions and the Lipschitz property of the systems in Equations \eqref{eq:model:cont} 
and \eqref{eq:adjointsys} (see References \cite{Jung:2002,SilvaTorres:2013} and references cited therein).
\end{Remark}

We implement Theorem~\ref{theo:PMP} numerically in Octave/MATLAB in
Section~\ref{sec:oc:numer}, and the optimal solution 
$\left(s^*(\cdot), i^*(\cdot), c^*(\cdot), a^*(\cdot)  \right)$ 
associated to the optimal control $u^*(\cdot)$ is computed 
for given parameter values and initial conditions.  


\subsection{Numerical Solution of the HIV Optimal Control Problem}
\label{sec:oc:numer}

The extremal given by Theorem~\ref{theo:PMP} is now computed numerically 
by implementing a forward-backward fourth-order Runge--Kutta method 
(see, e.g., Reference \cite{book:Lenhart}). This iterative method consists 
in solving the system in Equation \eqref{eq:model:cont} with a guess for the controls 
over the time interval $[0, T]$ using a forward fourth-order Runge--Kutta 
scheme and the transversality conditions $\lambda_i(T) = 0$, $i = 1, \ldots, 4$. 
Then, the adjoint system in Equation \eqref{eq:adjointsys} is solved by a backward 
fourth-order Runge--Kutta scheme using the current iteration solution 
of Equation \eqref{eq:model:cont}. The controls are updated by using a convex 
combination of the previous controls and the values from Equation \eqref{eq:opt:cont}. 
The iteration is stopped if the values of unknowns at the previous iteration 
are very close to the ones at the present iteration. Our 
\texttt{odeRungeKutta\_order4\_WithControl} function 
is implemented by the following GNU Octave instructions:


\medskip

{\small
	
\texttt{function dy = odeRungeKutta\_order4\_WithControl(T)}

\texttt{ \% Parameters of the model}

\texttt{ mi = 1.0 / 69.54; b = 2.1 * mi; beta = 1.6;}

\texttt{ etaC = 0.015; etaA = 1.3; fi = 1.0; ro = 0.1;}

\texttt{ alfa = 0.33; omega = 0.09; d = 1.0;}

\texttt{ }

\texttt{ \% Parameters of the Runge-Kutta (4th order) method}

\texttt{ test = -1; deltaError = 0.001; M = 1000;}

\texttt{ t = linspace(0,T,M+1);}

\texttt{ h = T / M; h2 = h / 2; h6 = h / 6;}

\texttt{ S = zeros(1,M+1); I = zeros(1,M+1);}

\texttt{ C = zeros(1,M+1); A = zeros(1,M+1);}

\texttt{ }

\texttt{ \% Initial conditions of the model}

\texttt{ S(1) = 0.6; I(1) = 0.2; C(1) = 0.1; A(1) = 0.1;}

\texttt{ }

\texttt{ \%Vectors for system restrictions and control}

\texttt{ Lambda1 = zeros(1,M+1); Lambda2 = zeros(1,M+1);}

\texttt{ Lambda3 = zeros(1,M+1); Lambda4 = zeros(1,M+1);}

\texttt{ U = zeros(1,M+1);}

\texttt{ \% Iterations of the method}

\texttt{ while(test $<$ 0)}

\texttt{ \ \ oldS = S; oldI = I; oldC = C; oldA = A;}

\texttt{ \ \ oldLambda1 = Lambda1; oldLambda2 = Lambda2;}

\texttt{ \ \ oldLambda3 = Lambda3; oldLambda4 = Lambda4;}

\texttt{ \ \ oldU = U;}

\texttt{ }

\texttt{ \ \ \%Forward Runge-Kutta iterations}

\texttt{ \ \ for i = 1:M}

\texttt{ \ \ \ \ \% Differential equations of the model}

\texttt{ \ \ \ \ \% First Runge-Kutta parameter}

\texttt{ \ \ \ \ aux1 = (1 - U(i)) * beta * (I(i) + etaC * C(i) + etaA * A(i)) * S(i);}

\texttt{ \ \ \ \ aux2 = d * A(i);}

\texttt{ }

\texttt{ \ \ \ \ auxS1 = b * (1 - S(i)) - aux1 + aux2 * S(i);}

\texttt{ \ \ \ \ auxI1 = aux1 - (ro + fi + b - aux2) * I(i) + alfa * A(i) + omega * C(i);}

\texttt{ \ \ \ \ auxC1 = fi * I(i) - (omega + b - aux2) * C(i);}

\texttt{ \ \ \ \ auxA1 = ro * I(i) - (alfa + b + d - aux2) * A(i);}

\texttt{ }

\texttt{ \ \ \ \ \% Second Runge-Kutta parameter}

\texttt{ \ \ \ \ auxU = 0.5 * (U(i) + U(i+1));}

\texttt{ \ \ \ \ auxS = S(i) + h2 * auxS1; auxI = I(i) + h2 * auxI1;}

\texttt{ \ \ \ \ auxC = C(i) + h2 * auxC1; auxA = A(i) + h2 * auxA1;}

\texttt{ \ \ \ \ aux1 = (1 - auxU) * beta * (auxI + etaC * auxC + etaA * auxA) * auxS;}

\texttt{ \ \ \ \ aux2 = d * auxA;}

\texttt{ \ \ \ \ auxS2 = b * (1 - auxS) - aux1 + aux2 * auxS;}

\texttt{ \ \ \ \ auxI2 = aux1 - (ro + fi + b - aux2) * auxI + alfa * auxA + omega * auxC;}

\texttt{ \ \ \ \ auxC2 = fi * auxI - (omega + b - aux2) * auxC;}

\texttt{ \ \ \ \ auxA2 = ro * auxI - (alfa + b + d - aux2) * auxA;}

\texttt{ }

\texttt{ \ \ \ \ \% Third Runge-Kutta parameter}

\texttt{ \ \ \ \ auxS = S(i) + h2 * auxS2; auxI = I(i) + h2 * auxI2;}

\texttt{ \ \ \ \ auxC = C(i) + h2 * auxC2; auxA = A(i) + h2 * auxA2;}

\texttt{ \ \ \ \ aux1 = (1 - auxU) * beta * (auxI + etaC * auxC + etaA * auxA) * auxS;}

\texttt{ \ \ \ \ aux2 = d * auxA;}

\texttt{ }

\texttt{ \ \ \ \ auxS3 = b * (1 - auxS) - aux1 + aux2 * auxS;}

\texttt{ \ \ \ \ auxI3 = aux1 - (ro + fi + b - aux2) * auxI + alfa * auxA + omega * auxC;}

\texttt{ \ \ \ \ auxC3 = fi * auxI - (omega + b - aux2) * auxC;}

\texttt{ \ \ \ \ auxA3 = ro * auxI - (alfa + b + d - aux2) * auxA;}

\texttt{ }

\texttt{ \ \ \ \ \% Fourth Runge-Kutta parameter}

\texttt{ \ \ \ \ auxS = S(i) + h * auxS3; auxI = I(i) + h * auxI3;}

\texttt{ \ \ \ \ auxC = C(i) + h * auxC3; auxA = A(i) + h * auxA3;}

\texttt{ \ \ \ \ aux1 = (1 - U(i+1)) * beta * (auxI + etaC * auxC + etaA * auxA) * auxS;}

\texttt{ \ \ \ \ aux2 = d * auxA;}

\texttt{ }

\texttt{ \ \ \ \ auxS4 = b * (1 - auxS) - aux1 + aux2 * auxS;}

\texttt{ \ \ \ \ auxI4 = aux1 - (ro + fi + b - aux2) * auxI + alfa * auxA + omega * auxC;}

\texttt{ \ \ \ \ auxC4 = fi * auxI - (omega + b - aux2) * auxC;}

\texttt{ \ \ \ \ auxA4 = ro * auxI - (alfa + b + d - aux2) * auxA;}

\texttt{ }

\texttt{ \ \ \ \ \% Runge-Kutta new approximation }

\texttt{ \ \ \ \ S(i+1) = S(i) + h6 * (auxS1 + 2 * (auxS2 + auxS3) + auxS4);}

\texttt{ \ \ \ \ I(i+1) = I(i) + h6 * (auxI1 + 2 * (auxI2 + auxI3) + auxI4);}

\texttt{ \ \ \ \ C(i+1) = C(i) + h6 * (auxC1 + 2 * (auxC2 + auxC3) + auxC4);}

\texttt{ \ \ \ \ A(i+1) = A(i) + h6 * (auxA1 + 2 * (auxA2 + auxA3) + auxA4);}

\texttt{ \ \ end}

\texttt{ }

\texttt{ \ \ \%Backward Runge-Kutta iterations}

\texttt{ \ \ for i = 1:M}

\texttt{ \ \ \ \ j = M + 2 - i;}

\texttt{ }

\texttt{ \ \ \ \ \% Differential equations of the model}

\texttt{ \ \ \ \ \% First Runge-Kutta parameter}

\texttt{ \ \ \ \ auxU = 1 - U(j);}

\texttt{ \ \ \ \ aux1 = auxU * beta * (I(j) + etaC * C(j) + etaA * A(j));}

\texttt{ \ \ \ \ aux2 = d * A(j);}

\texttt{ }

\texttt{ \ \ \ \ auxLambda11 = -1 + Lambda1(j) * (b + aux1 - aux2) - Lambda2(j) * aux1;}

\texttt{ \ \ \ \ aux1 = auxU * beta * S(j);}

\texttt{ \ \ \ \ auxLambda21 = 1 + Lambda1(j) * aux1 - Lambda2(j) * (aux1 - (ro + fi + b) + ...}

\texttt{ \ \ \ \ \ \ \ \ \ \ \ \ \ \ \ \ \ \ + aux2) - Lambda3(j) * fi -
	Lambda4(j) * ro;}

\texttt{ \ \ \ \ aux1 = auxU * beta * etaC * S(j);}

\texttt{ \ \ \ \ auxLambda31 = Lambda1(j) * aux1 - Lambda2(j) * (aux1 + ...}

\texttt{ \ \ \ \ \ \ \ \ \ \ \ \ \ \ \ \ \ \ + omega) + Lambda3(j) * (omega +
	b - aux2);}

\texttt{ \ \ \ \ aux1 = auxU * beta * etaA * S(j);}

\texttt{ \ \ \ \ auxLambda41 = Lambda1(j) * (aux1 + d * S(j)) ...}

\texttt{ \ \ \ \ \ \ \ \ \ \ \ \ \ \ \ \ \ \ - Lambda2(j) * (aux1 + alfa +
	...}

\texttt{ \ \ \ \ \ \ \ \ \ \ \ \ \ \ \ \ \ \ + d * I(j)) - Lambda3(j) * d *
	C(j) + ...}

\texttt{ \ \ \ \ \ \ \ \ \ \ \ \ \ \ \ \ \ \ + Lambda4(j) * (alfa + b + d - 2 * aux2);}

\texttt{ }

\texttt{ \ \ \ \ \% Second Runge-Kutta parameter}

\texttt{ \ \ \ \ auxU = 1 - 0.5 * (U(j) + U(j-1));}

\texttt{ \ \ \ \ auxS = 0.5 * (S(j) + S(j-1));}

\texttt{ \ \ \ \ auxI = 0.5 * (I(j) + I(j-1));}

\texttt{ \ \ \ \ auxC = 0.5 * (C(j) + C(j-1));}

\texttt{ \ \ \ \ auxA = 0.5 * (A(j) + A(j-1));}

\texttt{ }

\texttt{ \ \ \ \ aux1 = auxU * beta * (auxI + etaC * auxC + etaA * auxA);}

\texttt{ \ \ \ \ aux2 = d * auxA;}

\texttt{ \ \ \ \ auxLambda1 = Lambda1(j) - h2 * auxLambda11;}

\texttt{ \ \ \ \ auxLambda2 = Lambda2(j) - h2 * auxLambda21;}

\texttt{ \ \ \ \ auxLambda3 = Lambda3(j) - h2 * auxLambda31;}

\texttt{ \ \ \ \ auxLambda4 = Lambda4(j) - h2 * auxLambda41;}

\texttt{ }

\texttt{ \ \ \ \ auxLambda12 = -1 + auxLambda1 * (b + aux1 - aux2) - auxLambda2 * aux1;}

\texttt{ \ \ \ \ aux1 = auxU * beta * auxS;}

\texttt{ \ \ \ \ auxLambda22 = 1 + auxLambda1 * aux1 - auxLambda2 * (aux1 - (ro + fi + b) + ...}

\texttt{ \ \ \ \ \ \ \ \ \ \ \ \ \ \ \ \ \ \ + aux2) - auxLambda3 * fi - auxLambda4 * ro;}

\texttt{ \ \ \ \ aux1 = auxU * beta * etaC * auxS;}

\texttt{ \ \ \ \ auxLambda32 = auxLambda1 * aux1 - auxLambda2 * (aux1 + ...}

\texttt{ \ \ \ \ \ \ \ \ \ \ \ \ \ \ \ \ \ \ + omega) + auxLambda3 * (omega +
	b - aux2);}

\texttt{ \ \ \ \ aux1 = auxU * beta * etaA * auxS;}

\texttt{ \ \ \ \ auxLambda42 = auxLambda1 * (aux1 + d * auxS) ...}

\texttt{ \ \ \ \ \ \ \ \ \ \ \ \ \ \ \ \ \ \ - auxLambda2 * (aux1 + alfa +
	...}

\texttt{ \ \ \ \ \ \ \ \ \ \ \ \ \ \ \ \ \ \ + d * auxI) - auxLambda3 * d *
	auxC + ...}

\texttt{ \ \ \ \ \ \ \ \ \ \ \ \ \ \ \ \ \ \ + auxLambda4 * (alfa + b + d - 2 * aux2);}

\texttt{ }

\texttt{ \ \ \ \ \% Third Runge-Kutta parameter}

\texttt{ \ \ \ \ aux1 = auxU * beta * (auxI + etaC * auxC + etaA * auxA);}

\texttt{ \ \ \ \ auxLambda1 = Lambda1(j) - h2 * auxLambda12;}

\texttt{ \ \ \ \ auxLambda2 = Lambda2(j) - h2 * auxLambda22;}

\texttt{ \ \ \ \ auxLambda3 = Lambda3(j) - h2 * auxLambda32;}

\texttt{ \ \ \ \ auxLambda4 = Lambda4(j) - h2 * auxLambda42;}

\texttt{ }

\texttt{ \ \ \ \ auxLambda13 = -1 + auxLambda1 * (b + aux1 - aux2) - auxLambda2 * aux1;}

\texttt{ \ \ \ \ aux1 = auxU * beta * auxS;}

\texttt{ \ \ \ \ auxLambda23 = 1 + auxLambda1 * aux1 ...}

\texttt{ \ \ \ \ \ \ \ \ \ \ \ \ \ \ \ \ \ \ - auxLambda2 * (aux1 - (ro + fi +
	b) + ...}

\texttt{ \ \ \ \ \ \ \ \ \ \ \ \ \ \ \ \ \ \ + aux2) - auxLambda3 * fi -
	auxLambda4 * ro;}

\texttt{ \ \ \ \ aux1 = auxU * beta * etaC * auxS;}

\texttt{ \ \ \ \ auxLambda33 = auxLambda1 * aux1 - auxLambda2 * (aux1 + ...}

\texttt{ \ \ \ \ \ \ \ \ \ \ \ \ \ \ \ \ \ \ + omega) + auxLambda3 * (omega +
	b - aux2);}

\texttt{ \ \ \ \ aux1 = auxU * beta * etaA * auxS;}

\texttt{ \ \ \ \ auxLambda43 = auxLambda1 * (aux1 + d * auxS) ...}

\texttt{ \ \ \ \ \ \ \ \ \ \ \ \ - auxLambda2 * (aux1 + alfa  + d * auxI) - auxLambda3 * d * auxC + ...}

\texttt{ \ \ \ \ \ \ \ \ \ \ \ \ + auxLambda4 * (alfa + b + d - 2 * aux2);}

\texttt{ \ \ \ \ \% Fourth Runge-Kutta parameter}

\texttt{ \ \ \ \ auxU = 1 - U(j-1); auxS = S(j-1);}

\texttt{ \ \ \ \ auxI = I(j-1); auxC = C(j-1); auxA = A(j-1);}

\texttt{ }

\texttt{ \ \ \ \ aux1 = auxU * beta * (auxI + etaC * auxC + etaA * auxA);}

\texttt{ \ \ \ \ aux2 = d * auxA;}

\texttt{ \ \ \ \ auxLambda1 = Lambda1(j) - h * auxLambda13;}

\texttt{ \ \ \ \ auxLambda2 = Lambda2(j) - h * auxLambda23;}

\texttt{ \ \ \ \ auxLambda3 = Lambda3(j) - h * auxLambda33;}

\texttt{ \ \ \ \ auxLambda4 = Lambda4(j) - h * auxLambda43;}

\texttt{ }

\texttt{ \ \ \ \ auxLambda14 = -1 + auxLambda1 * (b + aux1 - aux2) ...}

\texttt{ \ \ \ \ \ \ \ \ \ \ \ \ \ \ \ \ \ \ - auxLambda2 * aux1;}

\texttt{ \ \ \ \ aux1 = auxU * beta * auxS;}

\texttt{ \ \ \ \ auxLambda24 = 1 + auxLambda1 * aux1 ...}

\texttt{ \ \ \ \ \ \ \ \ \ \ \ \ \ \ \ \ \ \ - auxLambda2 * (aux1 - (ro + fi +
	b) + ...}

\texttt{ \ \ \ \ \ \ \ \ \ \ \ \ \ \ \ \ \ \ + aux2) - auxLambda3 * fi -
	auxLambda4 * ro;}

\texttt{ \ \ \ \ aux1 = auxU * beta * etaC * auxS;}

\texttt{ \ \ \ \ auxLambda34 = auxLambda1 * aux1 - auxLambda2 * (aux1 + ...}

\texttt{ \ \ \ \ \ \ \ \ \ \ \ \ \ \ \ \ \ \ + omega) + auxLambda3 * (omega + b - aux2);}

\texttt{ \ \ \ \ aux1 = auxU * beta * etaA * auxS;}

\texttt{ \ \ \ \ auxLambda44 = auxLambda1 * (aux1 + d * auxS) ...}

\texttt{ \ \ \ \ \ \ \ \ \ \ \ \ \ \ \ \ \ \ - auxLambda2 * (aux1 + alfa +
	...}

\texttt{ \ \ \ \ \ \ \ \ \ \ \ \ \ \ \ \ \ \ + d * auxI) - auxLambda3 * d *
	auxC + ...}

\texttt{ \ \ \ \ \ \ \ \ \ \ \ \ \ \ \ \ \ \ + auxLambda4 * (alfa + b + d - 2
	* aux2);}

\texttt{ }

\texttt{ }

\texttt{ \ \ \ \ \% Runge-Kutta new approximation }

\texttt{ \ \ \ \ Lambda1(j-1) = Lambda1(j) - h6 * (auxLambda11 + ...}

\texttt{ \ \ \ \ \ \ \ \ \ \ \ \ \ \ \ \ \ \ \ + 2 * (auxLambda12 + auxLambda13) + auxLambda14);}

\texttt{ \ \ \ \ Lambda2(j-1) = Lambda2(j) - h6 * (auxLambda21 + ...}

\texttt{ \ \ \ \ \ \ \ \ \ \ \ \ \ \ \ \ \ \ \ + 2 * (auxLambda22 + auxLambda23) + auxLambda24);}

\texttt{ \ \ \ \ Lambda3(j-1) = Lambda3(j) - h6 * (auxLambda31 + ...}

\texttt{ \ \ \ \ \ \ \ \ \ \ \ \ \ \ \ \ \ \ \ + 2 * (auxLambda32 + auxLambda33) + auxLambda34); }

\texttt{ \ \ \ \ Lambda4(j-1) = Lambda4(j) - h6 * (auxLambda41 + ...}

\texttt{ \ \ \ \ \ \ \ \ \ \ \ \ \ \ \ \ \ \ \ + 2 * (auxLambda42 + auxLambda43) + auxLambda44);}

\texttt{ \ \ end}

\texttt{ }

\texttt{ \ \ \% New vector control}

\texttt{ \ \ for i = 1:M+1}

\texttt{ \ \ \ \ vAux(i) = 0.5 * beta * (I(i) + etaC * C(i) + ...}

\texttt{ \ \ \ \ \ \ \ \ \ \ \ \ \ \ + etaA * A(i)) * S(i) * (Lambda1(i) - Lambda2(i));}

\texttt{ \ \ \ \ auxU = min([max([0.0 vAux(i)]) 0.5]);}

\texttt{ \ \ \ \ U(i) = 0.5 * (auxU + oldU(i));}

\texttt{ \ \ end}

\texttt{ }

\texttt{ \ \ \% Absolute error for convergence}

\texttt{ \ \ temp1 = deltaError * sum(abs(S)) - sum(abs(oldS - S));}

\texttt{ \ \ temp2 = deltaError * sum(abs(I)) - sum(abs(oldI - I));}

\texttt{ \ \ temp3 = deltaError * sum(abs(C)) - sum(abs(oldC - C));}

\texttt{ \ \ temp4 = deltaError * sum(abs(A)) - sum(abs(oldA - A));}

\texttt{ \ \ temp5 = deltaError * sum(abs(U)) - sum(abs(oldU - U));}

\texttt{ \ \ temp6 = deltaError * sum(abs(Lambda1)) - sum(abs(oldLambda1 - Lambda1));}

\texttt{ \ \ temp7 = deltaError * sum(abs(Lambda2)) - sum(abs(oldLambda2 - Lambda2));}

\texttt{ \ \ temp8 = deltaError * sum(abs(Lambda3)) - sum(abs(oldLambda3 - Lambda3));}

\texttt{ \ \ temp9 = deltaError * sum(abs(Lambda4)) - sum(abs(oldLambda4 - Lambda4));}

\texttt{ \ \ test = min(temp1,min(temp2,min(temp3,min(temp4, ...}

\texttt{ \ \ \ \ \ \ \ \ \ min(temp5,min(temp6,min(temp7,min(temp8,temp9))))))));}

\texttt{ end}

\texttt{ dy(1,:) = t; dy(2,:) = S; dy(3,:) = I;}

\texttt{ dy(4,:) = C; dy(5,:) = A; dy(6,:) = U;}

\texttt{ }

\texttt{ disp("Value of LAMBDA at FINAL TIME");}

\texttt{ disp([Lambda1(M+1) Lambda2(M+1) Lambda3(M+1) Lambda4(M+1)]);}

}

\medskip

\begin{figure}[H]
\centering
\includegraphics[scale=0.38]{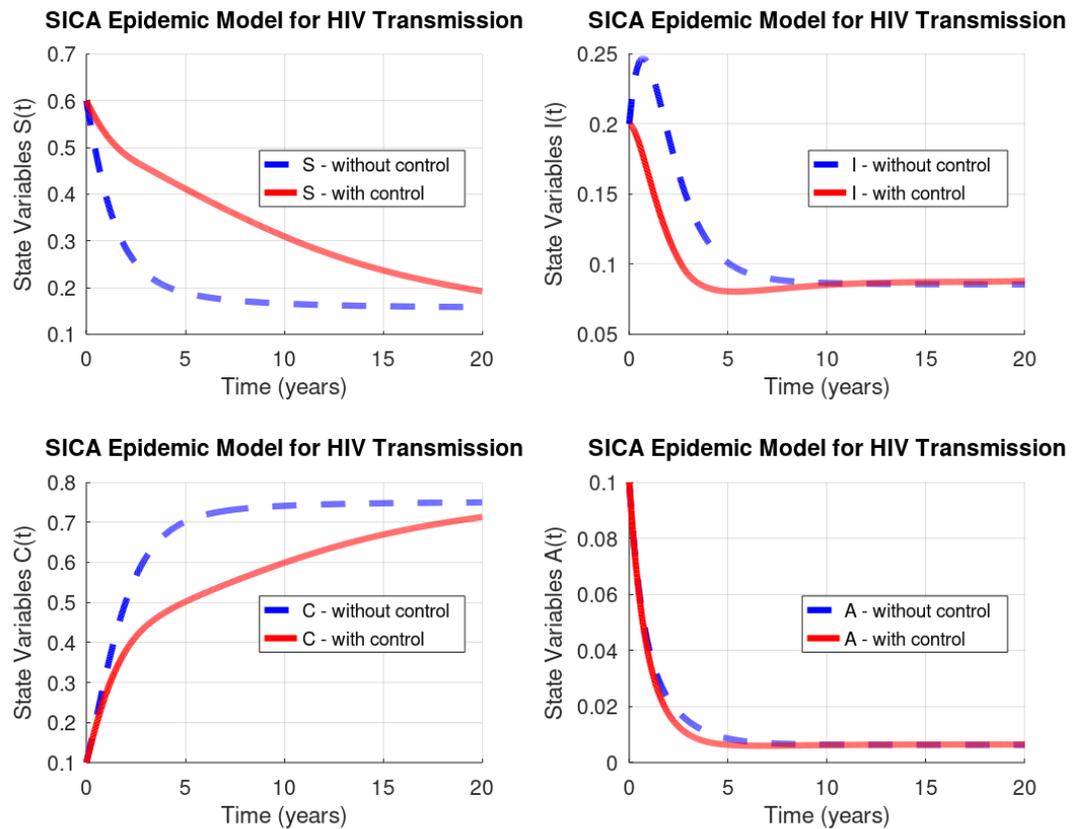}
\caption{Optimal state variables for the control problem in Equations
\eqref{eq:model:cont}--\eqref{eq:cost:max} subject to the initial 
conditions in Equation \eqref{initcond:num} and the admissible 
control set in Equation \eqref{eq:dmiss:cont:05}
\emph{versus} trajectories without control measures.}
\label{figura5}
\end{figure}
\begin{figure}[H]
\centering
\includegraphics[scale=0.25]{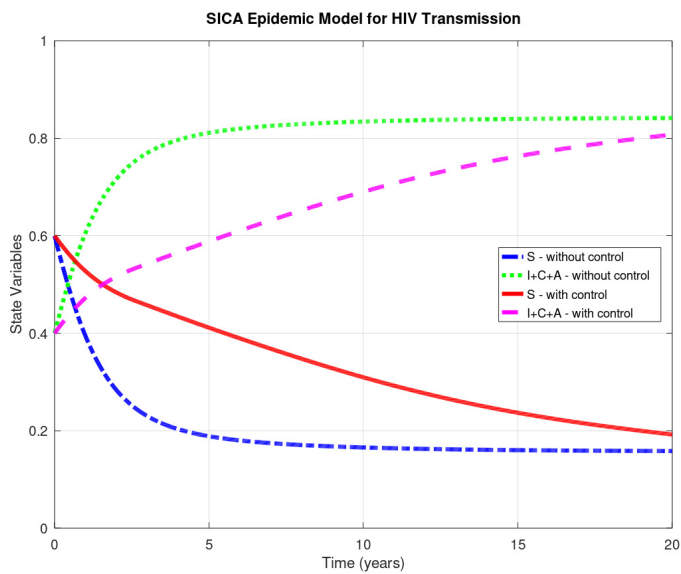}
\caption{Comparison: solutions to the initial value problem in Equations 
\eqref{eq:model:2}--\eqref{initcond:num} \emph{versus} solutions 
to the optimal control problem in Equations \eqref{eq:model:cont}--\eqref{eq:cost:max} 
subject to initial conditions in Equation \eqref{initcond:num} and the admissible control 
set in Equation \eqref{eq:dmiss:cont:05}.}
\label{figura6}
\end{figure}
\begin{figure}[H]
\centering
\includegraphics[scale=0.25]{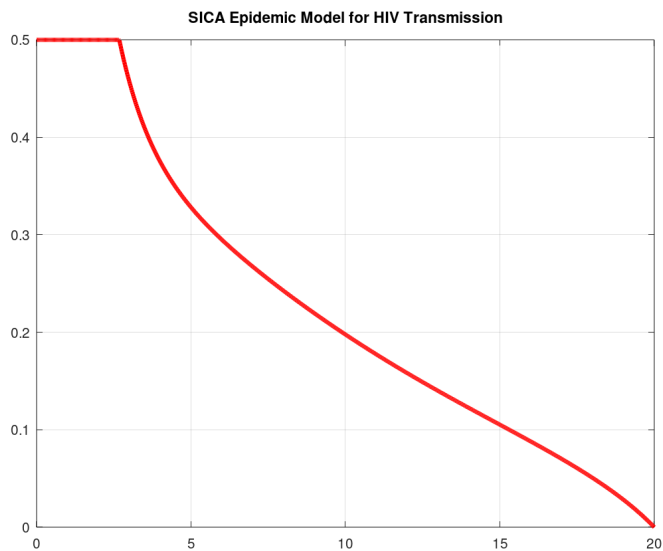}
\caption{Optimal control $u^*$ for the HIV optimal control problem in Equations 
\eqref{eq:model:cont}--\eqref{eq:cost:max} subject to the initial 
conditions in Equation \eqref{initcond:num} and the admissible control 
set in Equation \eqref{eq:dmiss:cont:05}.}
\label{figura7}
\end{figure}

For the numerical simulations, we consider $u_{max} = 0.5$, 
representing a lack of resources or misuse of the preventive HIV 
measures $u(\cdot)$, that is, the set of admissible controls is given by
\begin{equation}
\label{eq:dmiss:cont:05}
\Omega = \{ u(\cdot) \in L^\infty(0, T) \, | \, 0 
\leq u(t) \leq 0.5 \, , \, \forall t \in [0, T]  \} 
\end{equation}
with $T = 20$ (years). Figures~\ref{figura5}\ and \ref{figura6} show 
the numerical solution to the optimal control problem of Equations 
\eqref{eq:model:cont}--\eqref{eq:cost:max} with the initial 
conditions of Equation \eqref{initcond:num} and the admissible 
control set in Equation \eqref{eq:dmiss:cont:05} 
computed by our \texttt{odeRungeKutta\_order4\_WithControl} function. 
Figure~\ref{figura7}\ depicts the extremal control behaviour of $u^*$.


\section{Conclusion}

The paper provides a study on numerical methods 
to deal with modelling and optimal control of epidemic problems.
Simple but effective Octave/MATLAB code is fully provided
for a recent model proposed in Reference \cite{EcoComplexity2017}.
The given numerical procedures are robust with respect to the parameters: 
we have used the same values as the ones in Reference \cite{EcoComplexity2017},
but the code is valid for other values of the parameters
and easily modified to other models. The results show
the effectiveness of optimal control theory in medicine
and the usefulness of a scientific computing system such as GNU Octave:
using the control measure as predicted by Pontryagin's
maximum principle and numerically computed by our Octave code,
one sees that the number of HIV/AIDS-infected and -chronic individuals
diminish and, as a consequence, the number of susceptible (not ill)
individuals increase. We trust our paper will be very useful to a 
practitioner from the disease control area. Indeed, this work
has been motivated by many emails we received and continue to receive, 
asking us to provide the software code associated to our research papers 
on applications of optimal control theory in epidemiology, 
e.g., References \cite{MR3349757,MR3918295}.


\vspace{6pt} 

\funding{This research was partially supported by the Portuguese 
Foundation for Science and Technology (FCT) within projects 
UID/MAT/04106/2019 (CIDMA) and PTDC/EEI-AUT/2933/2014 (TOCCATTA) and 
was cofunded by FEDER funds through COMPETE2020---Programa Operacional Competitividade 
e Internacionaliza\c{c}\~ao (POCI) and by national funds (FCT). Silva is also supported 
by national funds (OE), through FCT, I.P., in the scope of the framework contract 
foreseen in numbers 4--6 of article 23 of the Decree-Law 57/2016 
of August 29, changed by Law 57/2017 of July 19.}


\acknowledgments{The authors are sincerely grateful to four anonymous 
reviewers for several useful comments, suggestions, and criticisms, which 
helped them to improve the paper.}


\authorcontributions{Each author equally contributed to this paper, 
read and approved the final manuscript.}


\conflictsofinterest{The authors declare no conflict of interest.} 


\reftitle{References}


\end{document}